\documentclass{amsart}
\usepackage[latin1]{inputenc}
\usepackage[french]{babel}
\usepackage{amsmath}
\usepackage{amscd}

\newcommand{\Sys}{\mathfrak{S}}
\newcommand{\R}{{\mathbb R}}
\newcommand{\C}{{\mathbb C}}

\newcommand{\vol}{{\rm vol}}
\newcommand{\sys}{{\rm sys}}

\newtheorem{main-theorem}{Theorem}
\newtheorem{theorem}{Th\'eor\`eme}

\newtheorem{proposition}{Proposition}
\newtheorem{question}{Question}

\theoremstyle{definition}

\title[Optimalit\'e systolique infinit\'esimale de l'oscillateur
harmonique]{Optimalit\'e systolique infinit\'esimale de l'oscillateur
  harmonique}

\author{J.C. \'Alvarez Paiva}

\address{J.C. \'Alvarez Paiva, Laboratoire Paul Painlev\'e, Bat. M2, Universit\'e des Sciences et Technologies,
59 655 Villeneuve d'Ascq, France.}

\email{juan-carlos.alvarez-paiva@math.univ-lille1.fr}

\author{F. Balacheff}

\address{F. Balacheff, Laboratoire Paul Painlev\'e, Bat. M2, Universit\'e des Sciences et Technologies,
59 655 Villeneuve d'Ascq, France.}

\email{florent.balacheff@math.univ-lille1.fr}



\begin{document}%

\begin{abstract}
Nous \'etudions les aspects infinit\'esimaux du probl\`eme suivant. 
Soit  $H$ un hamiltonien de $\R^{2n}$ dont la surface d'\'energie $\{H=1\}$ borde un domaine compact et \'etoil\'e de volume identique \`a celui de la boule unit\'e de $\R^{2n}$. La surface d'\'energie $\{H=1\}$ contient-elle une orbite p\'eriodique du syst\`eme hamiltonien 
$$
\left\{\begin{array}{ccc}
\dot{q}&= & \, \, \,  {\partial H \over \partial p} \\ 
\dot{p}&=& - {\partial H \over \partial q}
\end{array}\right.
$$
dont l'action soit au plus $\pi$ ?\\

\noindent \textsc{Abstract.} We study the infinitesimal aspects of the following problem. Let $H$ be a Hamiltonian on $\R^{2n}$ whose energy surface $\{H=1\}$ encloses a compact starshaped domain of volume equal to that of the unit ball in $\R^{2n}$. Does the energy surface $\{H=1\}$ carry a periodic orbit of the Hamiltonian system
$$
\left\{\begin{array}{ccc}
\dot{q}&= & \, \, \,  {\partial H \over \partial p} \\ 
\dot{p}&=& - {\partial H \over \partial q}
\end{array}\right.
$$
with action less than or equal to $\pi$ ?

\end{abstract}

\maketitle


Consid\'erons le hamiltonien homog\`ene de degr\'e $2$ de $\R^{2n}$ donn\'e par la formule 
$$
H_{st}(q_1,\ldots,q_n,p_1,\ldots,p_n)=\sum_{i=1}^n(q_i^2+p_i^2)
$$
et associ\'e au syst\`eme form\'e de deux oscillateurs harmoniques
d\'ecoupl\'es de m\^eme p\'eriode propre. Les orbites de ce syst\`eme
hamiltonien sont toutes p\'eriodiques de p\'eriode $\pi$. La surface
d'\'energie $\{H_{st}=1\}$ borde la boule unit\'e standard de $\R^{2n}$ dont le volume vaut $\pi^n\over n!$. Cet hamiltonien est un exemple d'hamiltonien {\it p\'eriodique}, {\it i.e.} un hamiltonien dont toutes les orbites sont p\'eriodiques et de m\^eme p\'eriode. Le type de question auquel nous nous int\'eressons ici est la suivante.

\begin{question}\label{q:1}  Si $H : \R^{2n} \to [0,\infty)$ d\'esigne un hamiltonien propre, lisse en dehors
  de l'origine et homog\`ene de degr\'e $2$, tel que le volume du domaine
  bord\'e par la surface d'\'energie $\{H= 1\}$ vaille $\pi^n \over n!$, existe-t-il une orbite p\'eriodique du syst\`eme hamiltonien 
$$
\left\{\begin{array}{ccc}
\dot{q}&= & \, \, \, {\partial H \over \partial p} \\
\dot{p}&=& - {\partial H \over \partial q}
\end{array}\right.
$$
dont la p\'eriode (co\"incidant avec l'action) soit  au plus $\pi$ ?
\end{question}

Il n'est pas clair que la r\'eponse \`a cette question soit positive en toute
g\'en\'eralit\'e, mais nous obtenons ici des r\'esultats allant dans ce sens pour un
hamiltonien suffisamment proche de l'hamiltonien $H_{st}$ (pour la topologie lisse). \\

Pour \'enoncer plus pr\'ecis\'ement nos r\'esultats, repla\c cons notre question dans un contexte plus g\'en\'eral. Soit $M^{2n-1}$ une vari\'et\'e de dimension impaire munie d'une forme de contact $\alpha$ ({\it i.e.} une $1$-forme $\alpha$ telle que $\alpha \wedge d\alpha^{n-1}$ soit une forme volume). Le {\it volume} de $(M,\alpha)$ est d\'efini comme la quantit\'e
$$
\vol(M,\alpha) = {1 \over n!} \int_M \alpha\wedge d\alpha^{n-1}.
$$
Rappelons qu'\'etant donn\'e une vari\'et\'e de contact, les orbites du champ de Reeb $X$ d\'efini par les \'equations
$$
d\alpha(X,\cdot)=0 \ \, \text{et} \ \, \alpha(X)=1
$$
s'appellent les {\it caract\'eristiques}. Nous d\'efinissons la {\it systole} par la formule
$$ 
\sys(M,\alpha)=\inf_\gamma \int_\gamma \alpha 
$$
o\`u l'infimum des actions est pris sur l'ensemble des caract\'eristiques ferm\'ees. S'il n'existe pas de caract\'eristiques ferm\'ees, nous poserons $\sys(M,\alpha)=\infty$. Enfin d\'efinissons le {\it volume systolique} comme la quantit\'e
$$
\Sys(M,\alpha)=\frac{\vol(M,\alpha)}{\sys(M,\alpha)^n}.
$$
Dans le cas  o\`u $(M,\alpha)$ est le cotangent unitaire d'une vari\'et\'e riemannienne muni de sa $1$-forme canonique, cette derni\`ere quantit\'e co\"incide \`a une constante dimensionnelle pr\`es avec la notion de volume systolique pr\'esent\'ee dans \cite{AlvarezBalacheff:09}. Remarquons que si $\phi : M \to N$ d\'esigne un diff\'eomorphisme entre deux vari\'et\'es de dimension impaire, et $\alpha$ une forme de contact sur $N$, alors $\Sys(M,\phi^\ast \alpha)=\Sys(N,\alpha)$.

D'apr\`es un r\'esultat de Taubes \cite{Tau07}, toute vari\'et\'e de contact
$(M,\alpha)$ de dimension $3$ admet une caract\'eristique ferm\'ee, et  donc
$\Sys(M,\alpha)>0$.  Dans le cas du cotangent unitaire d'une
vari\'et\'e de Finsler compacte, l'existence d'une caract\'eristique ferm\'ee d\'ecoule des travaux de Lusternik et Fet (voir \cite{Kli78}), et l'\'etude infinit\'esimale du volume systolique a \'et\'e initi\'ee dans \cite{AlvarezBalacheff:09}. Dans le cas d'un hypersurface compacte, lisse et \'etoil\'ee de $\R^{2n}$ (par \'etoil\'ee nous entendons que l'hypersurface borde
un domaine \'etoil\'e de $\R^{2n}$ contenant l'origine dans son int\'erieur), l'existence d'une caract\'eristique ferm\'ee est d\^ue \`a Rabinowitz \cite{Rab87}, et c'est cette derni\`ere situation qui nous int\'eresse plus particuli\`erement ici. \\

Nous travaillons dans l'espace euclidien standard $\R^{2n}$
muni des coordonn\'ees $(q_1,\ldots,q_n,p_1,\ldots,p_n)$, et nous noterons
$\alpha$ la $1$-forme ${1\over 2}\sum_{i=1}^n (p_idq_i-q_idp_i)$. Soit $\Sigma
\subset \R^{2n}$ une hypersurface compacte, lisse et \'etoil\'ee. La restriction de la
$1$-forme $\alpha$ \`a $\Sigma$ est une forme de contact. Nous d\'efinissons l'{\it hamiltonien associ\'e \`a $\Sigma$} comme l'unique fonction not\'ee $H_\Sigma : \R^{2n} \to [0,\infty)$ homog\`ene de degr\'e $2$ et lisse en dehors de l'origine  telle que $\Sigma$ co\"incide avec la surface de niveau $\{H_\Sigma= 1\}$. R\'eciproquement, \`a tout hamiltonien $H$ propre, lisse en dehors de l'origine et homog\`ene de degr\'e $2$, nous pouvons associer la surface de niveau $\{H=1\}$ qui est une hypersurface lisse et \'etoil\'ee. Remarquons que le hamiltonien associ\'e \`a la sph\`ere unit\'e standard $S^{2n-1}$ n'est autre que
$$
H_{st}=\sum_{i=1}^n(q_i^2+p_i^2).
$$
Les caract\'eristiques de $(\Sigma,\alpha_{|\Sigma})$ (les orbites du champ de
Reeb) co\"incident exactement avec les orbites du syst\`eme hamiltonien
correspondant \`a la surface d'\'energie $\{H=1\}$, et d'apr\`es  \cite{Rab87}, $(\Sigma,\alpha_{|\Sigma})$ admet au moins une caract\'eristique ferm\'ee. Dans ce contexte, nous pouvons
noter $\Sys(\Sigma)$ le volume systolique de $(\Sigma,\alpha_{|\Sigma})$. Remarquons que 
$$
\Sys(S^{2n-1})={1\over n!}.
$$
Nous reformulons la question \ref{q:1} de la mani\`ere suivante :
\begin{question}\label{q:2} Pour quelles hypersurfaces $\Sigma$ lisses \'etoil\'ees de $\R^{2n}$ l'in\'egalit\'e 
\begin{equation}
\label{eq}
\Sys(\Sigma)\geq {1\over n!}
\end{equation}
est-elle v\'erifi\'ee ?
\end{question}
Cette question constitue une variation autour d'un probl\`eme pos\'e par Viterbo dans \cite{Vit00}. En effet, dans le cas o\`u $\Sigma$ est convexe, minorer le volume systolique de $\Sigma$ \'equivaut \`a minorer le volume symplectique du domaine $K$ compact bordant $\Sigma$ par sa capacit\'e d'Ekland-Hofer-Zender (voir \cite{MS98}).  D'apr\`es \cite{AMO08}, si $\Sigma$ est convexe, il existe une constante $c>0$ ind\'ependante de la dimension telle que 
$$
\Sys(\Sigma)\geq {c\over n!}.
$$

Nous commen\c cons par montrer que toute hypersurface $\Sigma$ invariante par le flot du hamiltonien $H_{st}$ v\'erifie
$$
\Sys(\Sigma)\geq {1\over n!}.
$$
Ces hypersurfaces invariantes par le flot du hamiltonien $H_{st}$ correspondent aux bords des domaines $K$ {\it circulaires}, {\it i.e.} des domaines tels que $e^{i\theta} \, K=K$ pour tout r\'eel $\theta$ ($\R^{2n}$ \'etant identifi\'e $\C^n$). 
En particulier, les bords des domaines de Reinhardt (comme par exemple les domaines construits par Hermann dans \cite{Her98}) sont circulaires, et satisfont l'in\'egalit\'e systolique (\ref{eq}). Plus g\'en\'eralement, nous prouvons le r\'esultat suivant.

\begin{proposition}\label{prop1}
Soient $\Sigma$ et $\Sigma_0$ deux hypersurfaces lisses et \'etoil\'ees. Supposons que le flot de Reeb de $\Sigma_0$ soit p\'eriodique de p\'eriode $\sys(\Sigma_0)$ et que le flot du hamiltonien $H_\Sigma$ associ\'e \`a $\Sigma$ soit constant le long des orbites du flot hamiltonien $H_0$ associ\'e \`a $\Sigma_0$. Alors 
$$
\Sys(\Sigma)\geq \Sys(\Sigma_0).
$$
De plus, l'\'egalit\'e n'est possible que si $\Sigma$ et $\Sigma_0$ sont homoth\'etiques. 
\end{proposition}

\noindent {\it D\'emonstration.} Nous proc\'edons de mani\`ere analogue \`a la
preuve du Th\'eor\`eme 3.1 de \cite{AlvarezBalacheff:09}. Quitte \`a transformer $\Sigma$ par une homoth\'etie, nous pouvons supposer que
$$
\vol(\Sigma)=\vol(\Sigma_0).
$$
Soit  $\rho : \Sigma_0 \to (0,\infty)$ la fonction radiale d\'efinie par 
$$
\rho(q,p)={1\over \sqrt{H(q,p)}}
$$
 et soit $\delta : \Sigma_0 \to \Sigma$ l'application $(q,p) \mapsto \rho(q,p)(q,p)$. Nous noterons abusivement par $\alpha$ la restriction de $\alpha$ \`a $\Sigma_0$ et $\Sigma$. Nous avons $\delta^* \alpha = \rho \alpha$, d'o\`u 
$$
\vol(\Sigma) = \frac{1}{n!} \int_{\Sigma_0} \rho^n \alpha \wedge (d\alpha)^{n-1}.
$$
Comme $\vol(\Sigma)=\vol(\Sigma_0)$, la moyenne de $\rho^n$
sur $\Sigma_0$ est \'egale \`a $1$ et donc le minimum de $\rho$ est au plus $1$.

Soit $(q,p)$ un point de $\Sigma_0$ o\`u $\rho$ atteint son minimum et soit $\gamma$ la caract\'eristique ferm\'ee issue de ce point.
L'image de $\gamma$ par l'application radiale $\delta$ est une caract\'eristique ferm\'ee de $\Sigma$ et son action vaut $\min(\rho) \cdot  \sys(\Sigma_0) \leq \sys(\Sigma_0)$ (voir \cite{AlvarezBalacheff:09}). Le cas d'\'egalit\'e impose que $\rho=1$ et donc $\Sigma=\Sigma_0$.  
\qed

Dire que le hamiltonien $H$ est constant le long des orbites de $H_0$ est \'equivalent \`a dire que le crochet de Poisson
$$
\{H, H_0 \} := dH(X_{H_0})=\sum_{i=1}^n \left({\partial H \over \partial q_i}{\partial H_0 \over \partial p_i} - {\partial H \over \partial p_i}{\partial H_0 \over \partial q_i} \right)
$$
est identiquement nul. De m\^eme que dans \cite{AlvarezBalacheff:09}, nous pouvons
adapter de mani\`ere directe l'approche de Cushman dans \cite{Cushman:92} des
formes normales des syst\`emes hamiltoniens pour prouver le r\'esultat suivant :

\begin{proposition}\label{prop2}
Soit $\{H_t\}$ une famille lisse de hamiltoniens lisses en dehors de l'origine et homog\`enes de degr\'e $2$ telle que $H_0$ soit p\'eriodique. Pour chaque entier $N$, il existe une isotopie $\phi^{(N)}_t : \R^{2n}\setminus 0  \to \R^{2n}\setminus 0$ de l'identit\'e par des transformations symplectiques lisses et homog\`enes de degr\'e $1$ telle que  
$$
H_t \circ \phi^{(N)}_t= H_0 + tE_1 + \cdots + t^N E_N + o(t^N) ,
$$
o\`u $E_i : \R^{2n} \to \R$ sont des hamiltoniens lisses en dehors de l'origine et homog\`enes de degr\'e $2$ tels que $\{E_i,H_0\} = 0$ $(1 \leq i \leq N)$.
\end{proposition}


Nous pouvons alors montrer le th\'eor\`eme suivant, analogue du th\'eor\`eme 2.1 de \cite{AlvarezBalacheff:09} dans notre situation :
\begin{theorem}
Soit $\Sigma_0$ une hypersurface lisse et \'etoil\'ee dont le flot de Reeb  est p\'eriodique de p\'eriode $\sys(\Sigma_0)$.
Soit $N$ un entier quelconque. Fixons un $N$-jet de d\'eformation de $H_0$, soit $N$ fonctions $H_1,\ldots,H_N$ lisses en dehors de l'origine et homog\`enes de degr\'e $2$ quelconques. Alors il existe une d\'eformation lisse $H_t:=H_{\Sigma_t}$ de $H_0$ telle que
$$
H_t= H_0+tH_1+\ldots+t^NH_N+o(t^N)
$$
et
$$
\Sys(\Sigma_t)\geq \Sys(\Sigma_0).
$$
\end{theorem}

\noindent {\it D\'emonstration.} 
D'apr\`es la proposition \ref{prop2}, il existe une d\'eformation lisse $\phi^{(N)}_t : \R^{2n}\setminus 0  \to \R^{2n}\setminus 0$ de l'application identit\'e consistant en des transformations symplectiques lisses et homog\`enes de degr\'e $1$ telle que  
$$
(H_0+tH_1+\ldots+t^NH_N) \circ \phi^{(N)}_t= H_0 + tE_1 + \ldots + t^N E_N + o(t^N) ,
$$
avec $\{E_i,H_0\} = 0$ $(1 \leq i \leq N)$.

D\'efinissons une famille de hamiltoniens homog\`enes de degr\'e $2$ et lisses en dehors de l'origine par la formule
$$
H_t=(H_0 + tE_1 + \ldots + t^N E_N )\circ(\phi^{(N)}_t)^{-1}.
$$
Comme $\{H_t \circ \phi^{(N)}_t,H_0\}=0$, nous d\'eduisons de la proposition \ref{prop1} et du fait que $(\phi^{(N)}_t)^\ast \alpha = \alpha$ que
$$
\Sys(\Sigma_{H_t})= \Sys(\phi^{(N)}_t(\Sigma_{H_t}))= \Sys(\Sigma_{H_t \circ \phi^{(N)}_t})\geq \Sys(\Sigma_0).
$$
\qed



\end{document}